\documentclass[11pt]{amsart}

\usepackage{amsmath, amsfonts, amssymb,amsthm}
\usepackage{amstext}
\usepackage{mathrsfs}
\usepackage{float,epsf,subfigure}
\usepackage[all,cmtip]{xy}
\usepackage{hyperref}
\usepackage{algorithm}
\usepackage{algorithmic}
\usepackage{mathtools}
\usepackage{float}
\usepackage{bm}

\makeatletter
\@namedef{subjclassname@2020}{\textup{2020} Mathematics Subject Classification}
\makeatother

\newtheorem{theorem}{Theorem}[section]

\newtheorem{cor}[theorem]{Corollary}
\newtheorem{def-thm}[theorem]{Definition-Theorem}
\newtheorem{lemma}[theorem]{Lemma}

\newtheorem*{tha}{Theorem A}
\newtheorem*{thb}{Theorem B}

\theoremstyle{plain}
\theoremstyle{definition}

\newcommand{\C}{\mathbb C}

\begin{document}
\title[]{Entire solutions of a special binomial differential equation}
\author[ ]{Dong-hai Zhao, Jie Zhang and Liang-wen Liao}
\address{School of Mathematics\\ Nanjing University \\ Nanjing 210093, P. R. China}
\email{donghaizhao2021@163.com}
\address{School of Mathematics \\ China University of Mining and Technology \\ Xuzhou 221116, P. R. China}
\email{zhangjie1981@cumt.edu.cn}
\address{School of Mathematics \\ Nanjing University \\ Nanjing 210093, P. R. China}
\email{maliao@nju.edu.cn}
\date{}

\maketitle \thispagestyle{empty} \setcounter{page}{1}

\begin{abstract}
We consider the binomial differential equation
$$ ff'-a(f'')^2=be^{2c},$$
where $a(z)$, $b(z)$, $c(z)$ are polynomials such that $b(z)\not\equiv 0$ and $c(z)$ is non-constant. In this paper, we determine all entire solutions of this equation.
\end{abstract}

\vskip\baselineskip
\setlength\arraycolsep{2pt}
\medskip

\section{Introduction and main result}

Gundersen and Yang \cite{Gundersen-Yang} considered the following three types of nonlinear binomial differential equations:
$$ff''-a(f')^2=be^{2c},\quad ff'-a(f'')^{2}=be^{2c},\quad f'f''-af^{2}=be^{2c},$$ where $a, b, c$ are polynomials. For the first equation, they proved the following result.

\begin{tha} {\rm\cite{Gundersen-Yang}} Let $f$ be an entire solution of a differential equation of the form
\begin{equation}\label{MM.1}
ff''-a(f')^2=b(z)e^{2c(z)},
\end{equation}
where $a(z)$, $b(z)$, $c(z)$ are polynomials such that $b(z)\not\equiv 0$ and $c(z)$ is non-constant. Then we have one of the following situations:
\\(1) $b(z)$ is non-constant and $f(z)=p(z)e^{c(z)}$, where $p(z)$ is a polynomial.
\\(2) $a(z)\equiv \tau$ for some constant $\tau\not=1$, $b(z)\equiv\mu$ for some constant $\mu\not=0$, $c(z)=\lambda z$ for some constant $\lambda\not=0$, and
$f(z)=\alpha e^{\lambda z},$ where $\alpha$ is a constant satisfying $\alpha^2\lambda^2(1-\tau)=\mu$.\\
(3) $a(z)\equiv 1$, $b(z)\equiv\mu$ for some constant $\mu\not=0$, $c(z)=\lambda z^2+\gamma z$, and
$f(z)=\alpha e^{\lambda z^2+\gamma z},$
where $\alpha$, $\lambda$, $\gamma$ are constants such that $2\lambda \alpha^2=\mu$.\\
(4) $a(z)\equiv1$, $b(z)\equiv\mu$ for some constant $\mu\not=0$, $c(z)=\lambda z$ for some constant $\lambda\not=0$, and one of the following three cases occur:\\
$f(z)=(\alpha z + \beta)e^{\lambda z},$ where $\alpha$ and $\beta$ are constants such that $\alpha^2=-\mu$;\\
$f(z)=\alpha e^{2\lambda z}+\beta,$ where $\alpha$ and $\beta$ are constants satisfying $4\alpha\beta\lambda^2=\mu$;\\
$f(z)=\alpha e^{\gamma z}+\beta e^{(2\lambda-\gamma)z},$ where $\alpha$, $\beta$, $\gamma$ are constants satisfying
$\gamma(2\lambda-\gamma)\not=0$ and $4\alpha\beta(\lambda-\gamma)^2=\mu.$
\end{tha}
 Gundersen and Yang \cite{Gundersen-Yang} also studied the remaining two types of equations. These problems are motivated by classical results on the zeros and poles of differential polynomials. Hayman \cite{Haman,Haman2} proved that a meromorphic function $f$ such that $f$ and $f''$ have only finitely many zeros and poles must have the form $f=Re^P$, where $R$ is rational and $P$ is a polynomial. Langley \cite{Langley1} later weakened the assumption on the poles. For a zero-free meromorphic function $f=1/g$, where $g$ is entire, the identity
 $$f''=\frac{2(g')^2-gg''}{g^3}$$
  motivates the study of zero-free binomial differential polynomials. Mues \cite{Mues} classified the entire functions $g$ for which $gg''-a(g')^2$ is zero-free, where $a\ne1$ is a constant,
while Langley \cite{Lang2,Lang3} obtained related results for meromorphic functions under additional assumptions. Further related results can be found in \cite{Ber,Mues2}.\par
Recently, Wu, Liu and Hu \cite{Wu-Liu-Hu} and Xiang and Long \cite{Xiang-Long}
investigated the entire solutions of the third equation, while Gao
and Song \cite{Gao1} and Wang and Zhang \cite{Wang-Zhang} studied
those of the second equation. Exponential polynomial solutions
were also considered by Gao and Song \cite{Gao1}. Yang, Liao and Lu
\cite{Yang-Liao-Lu} and Long, Xia and Xiang \cite{Long} also
investigated the following more general equations and obtained
explicit forms of their entire solutions:
\begin{equation}\label{5.1}
\begin{aligned}
a(z)ff''-b(z)(f')^{2} &= c(z)e^{2d(z)},\\
a(z)ff'-b(z)(f'')^{2} &= c(z)e^{2d(z)},\\
a(z)f'f''- b(z)f^{2} &= c(z)e^{2d(z)},
\end{aligned}
\end{equation}
where $a(z)$, $b(z)$ and $c(z)$ are polynomials and $d(z)$ is a
non-constant entire function.
\par For the second equation
\begin{equation}\label{1.2}
 ff'-a(z)(f'')^2=b(z)e^{2c(z)},
\end{equation}
where $a(z)$, $b(z)$, $c(z)$ are polynomials such that $b(z)\not\equiv 0$ and $c(z)$ is non-constant, Gao and Song \cite{Gao1} give the following result:
\begin{thb}
Let $f$ be an entire solution of Eq.\eqref{1.2}. Moreover\\
(1)If $(b'+2bc')f'-bf''\equiv 0$. Then $a(z)\equiv\tau$ for some constant $\tau\neq 0$, $b(z)\equiv\mu$ for some constant $\mu\neq 0$, $c(z)=2\lambda z$ for some constant $\lambda\neq 0$, and we have
$$f(z)=\lambda_1e^{2\lambda z}+\lambda_2,$$
where
$$1-8\tau\lambda^3=0 \quad\text{and}\quad \lambda_1\lambda_2=\frac{\mu}{2\lambda}.$$
(2)If $(b'+2bc')f'-bf''\not\equiv 0$ and
$$m\left(r,\frac{f'}{(b'+2bc')f'-bf''}\right)=S(r,f).$$
Then, we have two possibilities:
\\(i)$b(z)$ is non-constant and we have
$$f(z)=R(z)e^{c(z)},$$
where $R(z)$ is a polynomial satisfying
$$R(R'+Rc') -a\left(R''+2R'c'+Rc''+R(c')^2\right)^2=b.$$
\\(ii)
$b(z)\equiv\mu$ for some constant $\mu\neq 0$,
$a(z)\equiv\tau$ for some constant $\tau\neq 0$,
$c(z)=\lambda z$ for some constant $\lambda\neq 0$,
and one of the following three cases occur:
$$f(z)=\alpha e^{\lambda z},$$
where $\alpha$ is a constant satisfying $\lambda\alpha^2(1-\tau\lambda^3)=\mu$;
$$f(z)=c_1e^{\lambda z}+c_2e^{-2\lambda z},$$
where $c_1,c_2$ are non-zero constants such that
$$1+8\tau\lambda^3=0, \qquad \frac{9}{8}\lambda(c_1)^2=\mu;$$
$$f(z)=c_1e^{t_1z}+c_2e^{t_2z},$$
where $t_1,t_2$ are two distinct constants such that
$$t_1+t_2=2\lambda\quad\text{and}\quad (t_1)^3=(t_2)^3=\frac{1}{\tau}.$$
$c_1,c_2$ are non-zero constants such that $6c_1c_2\lambda=\mu$.
\end{thb}

In this paper, we investigate the solutions of equation \eqref{1.2} satisfying
\begin{equation}\label{1.3}
m\left(r,\frac{f}{f'}\right)
\neq O\left(\log(rT(r,f))\right).
\end{equation}
Our main result is as follows.
\begin{theorem}\label{t.main}
Let $f$ be an entire solution of \eqref{1.2}. If $f$ satisfies \eqref{1.3}, then $a,b$ are constants and there are constants
$\lambda\in\C\setminus\{0\}$ and $c_0\in\C$ such that
$$ a(z)\equiv\frac{1}{8\lambda^3}, \quad and\quad c(z)=\lambda z+c_0,$$
and
\begin{equation}\label{e.solutionclass}
 f(z)=Ae^{2c(z)}+\frac{b}{2\lambda A}.
\end{equation}
\end{theorem}
Actually, from the proof of Theorem B in \cite{Gao1}, we found that the condition
$$m\left(r,\frac{f'}{(b'+2bc')f'-bf''}\right)=S(r,f)$$
implies that
$$m\left(r,\frac{f}{f'}\right)= O\left(\log(rT(r,f))\right).$$
Hence we obtain the following corollary:
\begin{cor}
All entire solutions of \eqref{1.2} are given by Theorem B.
\end{cor}

\par Throughout this paper, we may assume that readers are familiar with the standard notation and fundamental results in Nevanlinna theory, e.g., \cite{Ye,Hayman-book,Laine,Yang-Yi-book}. For example, a meromorphic function $a$ is called a small function of $f$ if and only if $T(r, a) = S(r, f )=o(T(r,f))$, $r\not\in E$, where $E$ is a set of finite measure, which may be different at each time. For a meromorphic function $f$, we denote by $\rho(f)$ the order of $f$.
\section{Some Lemmas}

\begin{lemma}\label{L.Laine}{\rm\cite{Laine}}
Let $f$ be a transcendental meromorphic function and $k\geq 1$ be an integer. Then
$$m\big(r,\frac{f^{(k)}}{f}\big)=O\big(\log T(r,f)+\log r\big)$$ possibly outside a set E of r with finite linear measure,
and if $f$ is of finite order of growth, then
$$m\big(r,\frac{f^{(k)}}{f}\big)=O(\log r).$$
\end{lemma}

\begin{lemma}\label{Clunie}{\rm\cite{Liao-Ye}}
Let $f$ be a transcendental meromorphic function in the complex plane and satisfying
$$f^nP(z,f)=Q(z,f),$$
where $P(z,f)$,$Q(z,f)$ are differential polynomials with rational function coefficients. If the degree of $Q(z,f)$ is at most $n$, then
$$m(r,P(z,f))=O(\log r)(r\to\infty),$$
if $f$ is of finite order, and
$$m(r,P(z,f))=O(\log (rT(r,f)))(r\to\infty),$$
possibly outside a set $E$ of $r$ with finite linear measure, if $f(z)$ is of infinite order.
\end{lemma}
 For an entire function $f$, denote by $\mu(r,f)$ the maximum term and $\nu(r,f)$ the central index of $f$.
\begin{lemma}\label{V-M1}{\rm\cite{Laine}}
If $g$ is an entire function of order $\rho$, then
$$\rho(g)=\limsup_{r\to\infty}\frac{\log^{+}\nu(r,g)}{\log r}=\limsup_{r\to\infty}\frac{\log^{+}\log^{+}\mu(r,g)}{\log r}.$$
\end{lemma}

\begin{lemma}\label{W-V}{\rm\cite{Laine}}
Let $g$ be a transcendental entire function, let $0<\delta<\frac{1}{4}$ and $z$ be such that $|z|=r$ and that
$$|g(z)|>M(r,g){\nu(r,g)}^{-\frac{1}{4}+\delta}$$
holds. Then there exists a set $F\subset \mathbb{R}_+$ of finite logarithmic measure, i.e. $\int_{F} dt/t < +\infty$, such that
$$g^{(m)}(z)=\left(\frac{\nu(r,g)}{z}\right)^m(1+o(1)) g(z)$$
holds for all $m\geq0$ and all $r\not\in F$.
\end{lemma}

\section{Proof of Theorem~\ref{t.main}}
\begin{proof}
Let $f$ be an entire solution of \eqref{1.2}.
\par We can assume that $a\not\equiv0$. In fact, if $a\equiv0$, \eqref{1.2} yields
$$ff'=be^{2c},\qquad (f^2)'=2be^{2c},$$
it follows from $b,c$ are polynomials that $\rho(f)<\infty$ and $f$ has only finitely many zeros. By the Hadamard factorization theorem,
$$f=Pe^Q,$$
where $P\not\equiv0$ and $Q$ are polynomials. Consequently,
$$\frac{f}{f'}=\frac{P}{P'+PQ'}$$
is a rational function, it follows that
$$m\left(r,\frac{f}{f'}\right)=O(\log r),$$
which is a contradiction to \eqref{1.3}.
We can claim that $f$ has order $\rho(f)=\deg c$.
\par In fact, let $k=\deg c\geq1$. Equation \eqref{1.2} implies
that $f$ is transcendental. If $\rho(f)<k$, then
$$\rho\bigl(ff'-a(f'')^2\bigr)<k,$$
which contradicts $\rho(be^{2c})=k$. Hence $\rho(f)\geq k$. \par
We next prove that $\rho(f)= k$.\\
Since equation \eqref{1.2} gives
$$\frac{f'}{f}-a\left(\frac{f''}{f}\right)^2=b\frac{e^{c}}{f^2},$$
Therefore Lemma \ref{W-V} yields that there is a set $E\subset[1,\infty)$ of finite logarithmic measure such that $|z_r|=r$, $|f(z_r)|=M(r,f)$ and
\begin{equation}\label{2.1}
\left(\frac{\nu(r,f)}{z_r}\right)(1+o(1))-a(z_r)\left(\frac{\nu(r,f)}{z_r}\right)^4(1+o(1))=b(z_r)\frac{e^{c(z_r)}}{(f(z_r))^2},
\end{equation}
as $r\to\infty$ with $r\notin E$. \par
If $\rho(f)>k$, then $e^c$ is a small function of $f$. Lemma \ref{V-M1} gives
$$\rho(f)=\limsup_{r\to\infty}\frac{\log\nu(r,f)}{\log r}>k,$$
which yields that there exists a sequence $r_n\to\infty$ such that and
$$\nu(r_n,f)>r_n^k.$$
Since $E$ has finite logarithmic measure, we can choose
$|z_r|\in[r_n,2r_n]\setminus E$ for all sufficiently large $n$.
By the monotonicity of $\nu$,
$$\frac{\nu(r)}{r}\geq\frac{\nu(r_n)}{2r_n}>\frac12r_n^{k-1}\to\infty.$$
It follows from \eqref{2.1} that the left-hand side tends to $\infty$ and right-hand side tends to 0, a contradiction. Hence we obtain $\rho(f)= k.$\par
The logarithmic derivatives of both sides of \eqref{1.2} give
$$\frac{(f')^2+ff''-a'(f'')^2-2af''f'''}{ff'-a(f'')^2}=\frac{b'+2bc'}{b},$$
which has the following form
$$f'(bf'-df)=f''[(a'b-ad)f''+2abf'''-bf],$$where
\begin{equation}\label{4.1}
d=b'+2bc'.
\end{equation}
Define a meromorphic function
\begin{equation}\label{4.2}
\xi=\frac{(a'b-ad)f''+2abf'''-bf}{f'}=\frac{bf'-df}{f''},
\end{equation}
it follows from $\rho(f)\leq k$ that $\rho(\xi)\leq k$.
The observation from \eqref{1.2} that $f'$ and $f''$ have at most finitely many common zeros results in that $\xi$ can only have finitely many poles at most, it follows that
\begin{equation}\label{4.1.1}
N(r,\xi) = O(\log r).
\end{equation}
From \eqref{4.2} we can obtain
\begin{equation}\label{4.2.1}
 \xi+b\frac{f}{f'} =(a'b-ad)\frac{f''}{f'}+2ab\frac{f'''}{f'}.
 \end{equation}
 Since $\rho(f)=k$, the Lemma \ref{L.Laine} yields that
$$ m\left(r,\frac{f''}{f'}\right)=O(\log r) \quad and \quad m\left(r,\frac{f'''}{f'}\right)=O(\log r).$$
It follows from \eqref{4.2.1} that
$$ m(r,\xi)=m\left(r,\frac{f}{f'}\right)+O(\log r).$$
Together with the estimate \eqref{4.1.1}, we obtain that
\begin{equation}\label{4.2.2}
 T(r,\xi)=m\left(r,\frac{f}{f'}\right)+O(\log r)
\end{equation}
Hence we can claim that $\xi$ is transcendental. In fact, if $\xi$ is rational, \eqref{4.2.2} gives $m\left(r,{f}/{f'}\right)=O(\log r)$ which is a contradiction to \eqref{1.3}.
\par Then \eqref{4.2} yields two following differential equations
\begin{equation}\label{4.3}
(a'b-ad)f''+2abf'''-bf-\xi f'=0
\end{equation}
and
\begin{equation}\label{4.4}
bf'-df-\xi f''=0.
\end{equation}
By differentiating \eqref{4.4}, we can obtain
\begin{equation}\label{4.5}
-\xi f'''+(b-\xi')f''+(b'-d)f'-d'f=0.
\end{equation}
 Eliminate $f'''$ from \eqref{4.3} and \eqref{4.5} yields
\begin{equation}\label{4.6}
Rf''+Sf'+b(2ad'+\xi)f=0,
\end{equation}
where
\begin{equation}\label{4.7}
R=2ab\xi'-2ab^2-a'b\xi+ad\xi
\end{equation}
and
\begin{equation}\label{4.8}
S=4ab^2c'+\xi^2.
\end{equation}
 Similarly eliminate $f''$ from \eqref{4.4} and \eqref{4.6} gives
\begin{equation}\label{4.9}
uf'=vf,
\end{equation}
where
\begin{equation}\label{4.10}
u=\xi S+bR
\end{equation}
and
\begin{equation}\label{4.11}
v=Rd-2a b d'\xi-b\xi^2.
\end{equation}
\par Differentiating \eqref{4.9} gives
\begin{equation}\label{4.9.1}
uf''=(v-u')f'+v'f.
\end{equation}
Eliminate $f'$ from \eqref{4.9} and \eqref{4.9.1} gives
\begin{equation}\label{4.9.2}
u^2f''=(v'u-vu'+v^2)f.
\end{equation}
Now substituting \eqref{4.9} and \eqref{4.9.2} into \eqref{4.4} yields
$$\left(\xi(v'u-vu'+v^2)-buv+du^2\right)f=0,$$
it follows from $f\not\equiv0$ that
\begin{equation}\label{4.12}
\xi(v'u-vu'+v^2)-buv+du^2=0.
\end{equation}
Substituting \eqref{4.7}-\eqref{4.8} into \eqref{4.10}-\eqref{4.11}, we obtain 
$$u=\xi^3+2ab^2\xi'+(4ab^2c'+abd-a'b^2)\xi-2ab^3,$$
and
$$v=-b\xi^2+2abd\xi'+(ad^2-a'bd-2abd')\xi-2ab^2d.$$
Then substituting this two equations into \eqref{4.12} yields
\begin{equation}\label{4.14}
 \xi^5b(\xi'+2c'\xi)=- Q(z,\xi,\xi',\xi''),
\end{equation}
where $Q$ is a differential polynomial of $\xi$ with polynomial coefficients and $\deg_\xi Q\leq 5$.
Put
\begin{equation}\label{4.13}
  H=\xi'+2c'\xi,
\end{equation}
we next prove that $H$ is rational. Since $\xi$ is transcendental, then the coefficients of \eqref{4.14} are small functions of $\xi$. From \eqref{4.14} and \eqref{4.13}, Lemma \ref{Clunie} yields that $m(r,bH)=O(\log r)$. Hence
$$ m(r,H)\leq m(r,bH)+m(r,1/b)=O(\log r).$$
Together with \eqref{4.1.1} and \eqref{4.13}, we obtain that
$$N(r,H)=O(\log r).$$
It follows that $T(r,H)=O(\log r)$, which shows that $H$ is rational.
\par Substituting \eqref{4.13} into \eqref{4.7} gives
\begin{equation}\label{4.15}
 R=2ab(H-b)+(ab'-a'b-2abc')\xi.
\end{equation}
By substituting \eqref{4.8} and \eqref{4.15} into \eqref{4.10}-\eqref{4.11}, we obtain that
\begin{equation}\label{4.16}
 u=\xi^3+L\xi+M,\qquad v=-b\xi^2+N\xi+P,
\end{equation}
where
$$L=2ab^2c'-a'b^2+abb', \quad M=2ab^2(H-b),$$
$$ N=(ab'-a'b-2abc')d-2abd', \quad P=2abd(H-b).$$
By differentiating \eqref{4.16}, \eqref{4.13} yields that
\begin{equation}\label{4.17}
 u'=-6c'\xi^3+3H\xi^2+(L'-2c'L)\xi+LH+M',
 \end{equation}
 and
 \begin{equation}\label{4.18}
 v'=(4bc'-b')\xi^2+(-2bH+N'-2c'N)\xi+NH+P'.
\end{equation}

From substituting \eqref{4.16}-\eqref{4.18} into \eqref{4.12}, it follows that
\begin{equation*}
\begin{aligned}
&\bigl[bH+N'+4c'N+2b^2\bigr]\xi^5\\
&+\Bigl[-2HN+bL'+(6bc'+b')L-3bN+P'+6c'P\Bigr]\xi^4\\
&+\Bigl[LN'-NL'+N^2+bM'+(8bc'+b')M+b^2L-bHL-3bP-3HP\Bigr]\xi^3\\
&+\Bigl[LP'-PL'+MN'-NM'+2c'(LP-MN)-2bHM+dL^2\\
&-bLN+b^2M+2NP\Bigr]\xi^2\\
&+\Bigl[MP'-PM'+H(MN-LP)+2dLM-b(LP+MN)+P^2\Bigr]\xi\\
&+dM^2-bMP=0,
\end{aligned}
\end{equation*}
it follows that
 \begin{equation}\label{4.19}
A_5\xi^5+A_4\xi^4+\cdots + A_0=0,
\end{equation}
where
 \begin{equation}\label{4.20}
A_5= bH+N'+4c'N+2b^2,
\end{equation}
 \begin{equation}\label{4.21}
A_4= -2HN+bL'+(6bc'+b')L-3bN+P'+6c'P=0,
\end{equation}
and $A_i$ $(i=0,\cdots,5)$ are rational.
\par Since $\xi$ is transcendental, by comparing the Nevanlinna characteristic functions the equation \eqref{4.19} yields that $A_i=0$ $(i=0,\cdots,5)$.

Next we show that $a$ and $c'$ are constant by comparing the leading terms. Put
$$m=\deg a,\qquad n=\deg b,\qquad s=\deg c'.$$
\par Suppose $m+3s>0$. From \eqref{4.1} and \eqref{4.16}, we obtain that
\begin{align}\label{4.22}
 d&=2bc'+O(z^{n+s-1}),\notag\\
 L&=2ab^2c'+O(z^{m+2n+s-1}),\\
 N&=-4ab^2(c')^2+O(z^{m+2n+2s-1}).\notag
 \end{align}
\par As $A_5=0$, \eqref{4.20} and \eqref{4.22} give
\begin{equation}\label{4.23}
 H=-\frac{N'+4c'N+2b^2}{b}
   =16ab(c')^3+O(z^{m+n+3s-1}).
\end{equation}
\par From \eqref{4.16} $P=2abd(H-b)$, it follows that
\begin{equation}\label{4.24}
P=4ab^2c'H+O(z^{2m+3n+4s-1}).
\end{equation}
Thus, \eqref{4.22} and \eqref{4.24} yield
 \begin{align}\label{4.25}
 -2HN+6c'P&=32ab^2(c')^2H+O(z^{2m+3n+5s-1}).\notag\\
 bL'&=O(z^{m+3n+s-1}),\\
 (6bc'+b')L-3bN&=O(z^{m+3n+2s}),\notag\\
 P'&=O(z^{2m+3n+4s-1}).\notag
 \end{align}
\par Consequently, $A_4=0$, \eqref{4.23} and \eqref{4.25} imply
$$ 0=A_4=512a^2b^3(c')^5+O(z^{2m+3n+5s-1}),$$
which is a contradiction to $m+3s>0$ and $b$ is nonzero. \\
\par Therefore $m=s=0$, which means that $a$ and $c'$ are nonzero constants. We take $\lambda$ and $c_0$ are constants such that
$$c'=\lambda\ne0, \quad and \quad c(z)=\lambda z+c_0.$$
Hence \eqref{4.16} yields
\begin{equation}\label{4.27}
 L=abd,\qquad N=a(b'-2\lambda b)d-2abd',\qquad P=2abd(H-b).
 \end{equation}
\par We next prove that $b$ is also a nonzero constant. From $A_5=0$, \eqref{4.20} and \eqref{4.27} yield
\begin{equation}\label{4.26}
 H=2(8a\lambda^3-1)b+24a\lambda^2b'+12a\lambda b''+2ab'''.
\end{equation}
Then $A_4=0$, \eqref{4.21} and \eqref{4.27} yields
\begin{equation}\label{4.28}
 0=A_4=ab\left(2dH'+(16\lambda d+6d')H+5(bd'-b'd)\right).
\end{equation}
Since $a$, $b$ are nonzero, \eqref{4.28} gives
\begin{equation}\label{4.29}
 2dH'+(16\lambda d+6d')H+5(bd'-b'd)=0.
\end{equation}
If $n>0$, we can claim that $8a\lambda^3-1=0$. If not, $8a\lambda^3-1\not=0$, then from \eqref{4.26} we $\deg H=\deg b=n>0$. Therefore equation \eqref{4.29} yields
$$32\lambda^2bH+O(z^{2n-1})=0,$$
 a contradiction.
 \par Hence \eqref{4.26} yields
 \begin{equation}\label{4.30}
 H=24a\lambda^2b'+12a\lambda b''+2ab'''\not\equiv0,
\end{equation}
where $\deg H=\deg b-1=n-1$.
\par By substituting \eqref{4.30} into \eqref{4.29}, we obtain that
$$32\lambda^2bH+O(z^{2n-2})=0,$$
a contradiction. Hence $n=0$, which means that $b$ is constant, then \eqref{4.26} yields $H$ is constant. Since $d$, $b$, $H$ are constants, \eqref{4.29} gives
$$16\lambda dH=0,$$
it follows that $H=0$ from $b,d\not=0$. Thus \eqref{4.26} gives that $a=1/{8\lambda^3}$.
\par Now \eqref{4.13} and $L,M,N,P$ reduce to
\begin{equation}\label{4.31}
 \xi'=-2\lambda\xi,
\end{equation}
and
\begin{equation}\label{4.31.1}
L=2ab^2\lambda, \quad M=-2ab^3, \quad N=-4ab^2\lambda^2, \quad P=-4ab^3\lambda.
\end{equation}
By substituting \eqref{4.31.1} and $8a\lambda^3-1=0$ into \eqref{4.16}, we obtain
\begin{equation}\label{4.33}
 u=\xi^3+\frac{b^2}{4\lambda^2}\xi-\frac{b^3}{4\lambda^3},
 \end{equation}
and
\begin{equation}\label{4.34}
 v=-b\xi^2-\frac{b^2}{2\lambda}\xi-\frac{b^3}{2\lambda^2},
 \end{equation}
it follows that
\begin{equation}\label{4.35}
 u=\left(\frac1{2\lambda}-\frac{\xi}{b}\right)v.
  \end{equation}
Since $\xi$ is transcendental, \eqref{4.33} and \eqref{4.34} lead to that $u$ and $v$ are also transcendental. Eliminate $u$, $v$ from \eqref{4.9} and \eqref{4.35} gives
\begin{equation}\label{4.32}
 (b-2\lambda\xi)f'=2\lambda bf.
\end{equation}
By differentiating \eqref{4.32} and using \eqref{4.31}, we obtain
$$(b-2\lambda\xi)(f''-2\lambda f')=0.$$
Since $b-2\lambda\xi$ is transcendental, it follows that
$$f''=2\lambda f',$$
and therefore
$$f=De^{2\lambda z}+E,$$
where $D\ne0$ and $E$ are constants.
Substituting $f$ into \eqref{1.2} gives
$$ff'-\frac{(f'')^2}{8\lambda^3}=2\lambda DEe^{2\lambda z}=be^{2c_0}e^{2\lambda z}.$$
Thus $2\lambda DE=be^{2c_0}$. Setting
$A=De^{-2c_0}\ne0$ yields
\begin{equation}\label{rewrite:e.finalsolution2}
 f(z)=Ae^{2c(z)}+\frac{b}{2\lambda A},\qquad A\ne0.
\end{equation}
This proves the theorem.
\end{proof}

\enlargethispage{4\baselineskip}

\section*{Declarations}
\noindent\textbf{Funding.}
This research is supported by the National Natural Science
Foundation of China (Grant Nos.~12471077 and 12171232).

\noindent\textbf{Conflict of interest.}
The authors declare that they have no conflict of interest.

\end{document}